\documentclass[10pt,reqno,a4paper]{amsart}

\usepackage[T1]{fontenc}
\usepackage[utf8]{inputenc}
\usepackage[english]{babel}
\usepackage[a4paper,margin=1.05in]{geometry}
\usepackage{microtype}
\usepackage{amsmath,amssymb,amsthm,mathtools,mathrsfs}
\usepackage[hidelinks]{hyperref}
\hypersetup{
  pdftitle={On Residue Localization for Complex Supermanifolds},
  pdfauthor={Leonardo Abath, Mauricio Correa, and Miguel Rodriguez Pena},
  pdfkeywords={complex supermanifold, odd holomorphic vector field, holomorphic localisation, maximal-picture integral form, Berezinian density, holomorphic residue}
}
\numberwithin{equation}{section}

\theoremstyle{plain}
\newtheorem{theorem}{Theorem}[section]
\newtheorem{proposition}[theorem]{Proposition}
\newtheorem{lemma}[theorem]{Lemma}
\newtheorem{corollary}[theorem]{Corollary}
\newtheorem*{maintheorem}{Main theorem}

\theoremstyle{definition}
\newtheorem{definition}[theorem]{Definition}
\newtheorem{convention}[theorem]{Convention}

\theoremstyle{remark}
\newtheorem{remark}[theorem]{Remark}

\newcommand{\cO}{\mathcal O}

\newcommand{\cJ}{\mathcal J}

\newcommand{\Res}{\operatorname{Res}}
\newcommand{\C}{\mathbb C}
\newcommand{\PP}{\mathbb P}
\newcommand{\red}{\mathrm{red}}
\newcommand{\supp}{\operatorname{supp}}
\newcommand{\Der}{\operatorname{Der}}
\newcommand{\cE}{\mathcal E}
\newcommand{\cT}{\mathcal T}

\newcommand{\ddbar}{\bar\partial}

\newcommand{\nil}{\mathrm{nil}}

\begin{document}

\title[Residue localization on supermanifolds]{On Residue Localization for Complex Supermanifolds}

\author[L. Abath]{Leonardo Abath}
\author[M. Corr\^ea]{Maur\'icio Corr\^ea}
\author[M. Rodr\'{i}guez]{Miguel Rodr\'iguez Pe\~na}

\address[L. Abath]{CIMAT, Ap. Postal 402, Guanajuato 36000, Gto., M\'exico}
\email{leoabath@gmail.com}

\address[M. Corr\^ea]{ 
Universit\`a degli Studi di Bari, 
Via E. Orabona 4, I-70125, Bari, Italy
}
\email[M. Corr\^ea]{mauricio.barros@uniba.it, mauriciomatufmg@gmail.com } 

\address[M. Rodr\'{i}guez]{DEMAT--UFSJ, Departamento de Matem\'atica e Estat\'istica, S\~ao Jo\~ao del-Rei MG, Brazil}
\email{miguel.rodriguez.mat@ufsj.edu.br}

\begin{abstract}
We prove a holomorphic residue localization formula for odd holomorphic vector fields on compact complex supermanifolds whose fermionic and bosonic dimensions coincide.  Under isolated non-degeneracy hypotheses on the reduced zero set,    we  give an explicit local residue formula. 
\end{abstract}

\subjclass[2020]{Primary 32C11, 58A50,32A27, 14M30, 32S65, 58A10}

\keywords{complex supermanifold, odd holomorphic vector field, holomorphic localisation, maximal-picture integral form, Berezinian density, holomorphic residue}

\maketitle

\section{Introduction}

Localisation formulae convert global integrals into local contributions supported at the zero set of a symmetry or vector field.  In the complex and equivariant settings this idea underlies the Duistermaat--Heckman theorem, the Atiyah--Bott and Berline--Vergne formulae, and Bott's residue theorem for holomorphic vector fields \cite{DH,AB,BV,Bott}.  It is also closely connected with the residue theories of Baum--Bott, Carrell--Lieberman and Liu \cite{BaumBott,CarrellLieberman,Liu}.  In supergeometry, odd vector fields provide a natural localisation mechanism; this point of view appears in Witten's supersymmetric interpretation of localisation and in the theorem of Schwarz and Zaboronsky \cite{WittenMorse,WittenNotes,SchwarzZaboronsky}.

We work throughout with complex supermanifolds in the usual sheaf-theoretic sense.  Thus odd differentiation is differentiation in the structure sheaf, not differentiation of non-canonical Grassmann-valued representatives.  The forms to be integrated are smooth integral forms of maximal picture, or equivalently local representatives of Berezinian densities.  The Stokes theorem and the Berezin trace used below are the standard ones for integral forms \cite{Berezin,BernsteinLeites,Voronov,WittenNotes}.

The localisation operator is the Dolbeault--Cartan operator attached to an odd holomorphic vector field.  The use of the Dolbeault operator is essential: adding contraction to the holomorphic de Rham operator gives the holomorphic Lie derivative.  In the present integral-form setting, however, one must not treat the operator as automatically square-zero on the whole algebra of integral forms.  The deformation argument below uses only the precise identity needed for localising one-forms, namely that \(D_V(D_V\omega)=0\).

The principal result is the following residue formula.  A reduced zero is called triangular non-degenerate if, in suitable local supercoordinates centred at the zero, the reduced vertical component of the vector field is represented by holomorphic functions of the reduced variables alone and has invertible Jacobian.  This hypothesis is local and is used only for the explicit point-residue computation.

\begin{maintheorem} \label{thm:introduction-main}
Let \(S\) be a compact complex supermanifold of dimension \(n|n\), and let \(V\) be an odd holomorphic vector field whose reduced zero locus is finite.  Assume that every reduced zero is triangular non-degenerate.  Then every smooth maximal-picture integral form \(\eta\) satisfying \((\bar\partial+i_V)\eta=0\) satisfies
\begin{equation}\label{eq:intro-main-formula}
\int_S\eta=
\left(\frac{2\pi}{i}\right)^n
\sum_{p\in Z(V)_{\red}}
\frac{\mathcal N_p(V,\eta)}{\det J_g(p)} .
\end{equation}
Here \(J_g(p)\) is the Jacobian of the reduced vertical equations in a triangular chart at \(p\), and \(\mathcal N_p(V,\eta)\) is the local Wick numerator defined in Definition~\ref{def:local-numerator-expanded}.
\end{maintheorem}

The numerator \(\mathcal N_p(V,\eta)\) is obtained by normal-ordering the local integral-form expansion, projecting to the fixed Berezinian symbol, and applying the normalised Wick functional for the Gaussian \(e^{-|g|^2/t}\).  In the basic case \(\eta=\Phi_p\Delta_p\), provided the representative is \(A_{\nil}\)-silent in the sense made precise below, it reduces to
\begin{equation}\label{eq:intro-basic-formula}
[\xi_1\cdots\xi_n\bar\xi_1\cdots\bar\xi_n]\Phi_p(p,\xi,\bar\xi).
\end{equation}
In general the terms involving \(A_{\nil}\) are not discarded; they are part of \(\mathcal N_p(V,\eta)\).

Formula \eqref{eq:intro-main-formula} is proved in two stages.  First, a Dolbeault--Cartan deformation argument gives a residue decomposition for arbitrary isolated reduced zeros.  Secondly, at a triangular non-degenerate zero, a local calculation separates the ordinary complex Gaussian integral from the finite algebra of integral forms.  The Gaussian integral produces the holomorphic determinant in the denominator, while the odd and integral-form variables produce the numerator \(\mathcal N_p(V,\eta)\).  When no derivative delta symbols occur and the nilpotent scalar part of \(i_V\omega\) is silent in the normalised Gaussian limit, the numerator is the Berezin coefficient displayed in \eqref{eq:intro-basic-formula}; otherwise the ordered projection is followed by the Wick functional, which keeps precisely the finite Gaussian moments produced by \(A_{\nil}\), by derivative delta symbols, and by the Taylor expansion of the holomorphic Jacobian factor.

The same deformation argument gives a vanishing theorem: if the reduced zero locus of \(V\) is empty, the integral of every \((\bar\partial+i_V)\)-closed maximal-picture integral form is zero.  The balanced dimension hypothesis is used precisely in the local computation of point residues.  In dimension \(n|n\), the reduced vertical component gives \(n\) holomorphic equations on the \(n\)-dimensional reduced manifold, so a non-degenerate zero has a square Jacobian and an isolated residue.  In unbalanced dimensions one expects positive-dimensional or excess residues rather than the point formula proved here.

Section~\ref{sec:conventions} fixes the conventions on supermanifolds and integral forms.  Section~\ref{sec:differential} introduces the Dolbeault--Cartan operator and localising one-forms.  Section~\ref{sec:global} proves the deformation argument, the vanishing theorem and the global residue decomposition.  Sections~\ref{sec:local}--\ref{sec:local-residue-formula} contain the local calculation.  Section~\ref{sec:explicit-local-models} records low-dimensional models, and Section~\ref{sec:applications} gives two applications.

\section{Conventions on supermanifolds and integral forms}\label{sec:conventions}

We use the standard definition of a complex supermanifold as a locally ringed space \cite{Leites,Manin,Varadarajan,Tuynman}.  Thus a complex supermanifold of dimension \(m|q\) is a pair \(S=(S_{\red},\cO_S)\) locally isomorphic to \((U,\cO_U\otimes\wedge^\bullet\C^q)\), with \(U\subset\C^m\) open.  The nilpotent ideal generated by the odd functions is denoted by \(\cJ_S\), so that \(\cO_{S_{\red}}=\cO_S/\cJ_S\).  The quotient \(\cE^\vee=\cJ_S/\cJ_S^2\) is the odd conormal bundle along the reduced manifold, and \(\cE\) denotes its dual.

In local supercoordinates \((z_1,\ldots,z_m;\xi_1,\ldots,\xi_q)\), every holomorphic superfunction has a unique finite expansion
\begin{equation}\label{eq:superfunction-expansion}
        h=\sum_I h_I(z)\xi_I,\qquad h_I\in\cO_U,
\end{equation}
where \(I\) runs through ordered subsets of \(\{1,\ldots,q\}\).  The odd derivative \(\partial/\partial\xi_j\) is the odd derivation of the local superalgebra determined by \(\partial\xi_k/\partial\xi_j=\delta_{jk}\) and by the super-Leibniz rule.  Thus all odd derivatives used below are derivations of \(\cO_S\).  No derivative is taken with respect to an ambiguous Grassmann-valued representative of a function.

The holomorphic tangent sheaf is \(\cT_S=\Der_\C(\cO_S)\).  An odd holomorphic vector field is a section \(V\in H^0(S,\cT_S)_{\bar1}\), and locally it has the form
\begin{equation}\label{eq:local-vector}
V=\sum_{i=1}^m f_i(z,\xi)\frac{\partial}{\partial z_i}
  +\sum_{\alpha=1}^q g_\alpha(z,\xi)\frac{\partial}{\partial\xi_\alpha},
\end{equation}
where the \(f_i\) are odd and the \(g_\alpha\) are even.  Reducing the vertical coefficients modulo \(\cJ_S\) gives a canonical section
\begin{equation}\label{eq:reduced-vertical-section}
        s_V\in H^0(S_{\red},\cE).
\end{equation}
The reduced zero locus \(Z(V)_{\red}\) is the zero scheme of this section.  In the balanced case \(m=q=n\), a reduced zero is called triangular non-degenerate if, in coordinates centred at the zero, the section \(s_V\) is represented by holomorphic functions \(g_1(z),\ldots,g_n(z)\) and the matrix \((\partial g_i/\partial z_j)\) is invertible at the zero.

We integrate smooth integral forms of maximal picture, equivalently smooth representatives of Berezinian densities \cite{Berezin,Manin,Varadarajan,BernsteinLeites,Voronov,WittenNotes}.  In local coordinates of dimension \(n|n\), set
\begin{equation}\label{eq:Delta-basic}
\Delta=\delta(d\xi_1)\cdots\delta(d\xi_n)
        \delta(d\bar\xi_1)\cdots\delta(d\bar\xi_n)
\end{equation}
and
\begin{equation}\label{eq:basic-integral-form}
\mathcal D=dz_1\wedge\cdots\wedge dz_n\wedge d\bar z_1\wedge\cdots\wedge d\bar z_n\,\Delta .
\end{equation}
For a local representative \(\Phi\mathcal D\), the integral is
\begin{equation}\label{eq:berezin-integral}
\int \Phi\mathcal D=
\int_{S_{\red}}[\xi_1\cdots\xi_n\bar\xi_1\cdots\bar\xi_n]\Phi\,
 dz_1\wedge\cdots\wedge dz_n\wedge d\bar z_1\wedge\cdots\wedge d\bar z_n.
\end{equation}
If an integral form is inhomogeneous in the ordinary differential degree, its integral means the integral of its component proportional to \(\mathcal D\).  This is the convention used in the finite expansions appearing in the localisation argument.

General maximal-picture representatives may involve derivatives of the delta symbols.  The local calculation uses only the standard identities
\begin{equation}\label{eq:integral-form-identities}
 d\xi_i\delta(d\xi_i)=0,\qquad
 d\xi_i\delta^{(r)}(d\xi_i)=-r\delta^{(r-1)}(d\xi_i),
\end{equation}
and the analogous formulae for \(d\bar\xi_i\).  The ordering in \eqref{eq:basic-integral-form}, together with the ordinary orientation convention implicit in \eqref{eq:berezin-integral}, fixes all signs in the coefficient extractions below.

\begin{convention}[Normal ordering and signs]\label{conv:normal-ordering}
All local products of integral forms are read in the ordered superalgebra of integral forms.  The coordinates \(z_i,\bar z_i\) are even, the coordinates \(\xi_i,\bar\xi_i\) are odd, and their differentials have the usual opposite parities.  In particular, \(dz_i,d\bar z_i\) are ordinary odd one-forms, whereas \(d\xi_i,d\bar\xi_i\) are even variables and are not exterior-nilpotent.  Thus no calculation below uses an identity such as \((d\xi_i)^2=0\) or \(C^2=0\).

To compute the coefficient of \(\mathcal D\), one first moves the ordinary differentials into the displayed order
\[
dz_1\cdots dz_n\,d\bar z_1\cdots d\bar z_n
\]
using the Koszul sign rule, then applies the distributional identities \eqref{eq:integral-form-identities} to the even variables \(d\xi_i,d\bar\xi_i\), and finally moves the surviving delta symbols into the order fixed in \eqref{eq:Delta-basic}.  Whenever a closed sign is not displayed, the sign is by definition the one obtained by this normal-ordering procedure.  No factor is commuted as if the algebra were an ordinary commutative algebra.
\end{convention}

\section{The localisation operator}\label{sec:differential}

Let \(V\) be an odd holomorphic vector field.  On smooth integral forms we use the Dolbeault--Cartan operator
\begin{equation}\label{eq:DV-def}
D_V=\ddbar+i_V.
\end{equation}
Here \(i_V\) denotes contraction by the superderivation \(V\).  We shall call an integral form \(D_V\)-closed if it is annihilated by \(D_V\).  This is a coordinate-free condition, because both \(\ddbar\) and contraction by the global vector field act naturally on integral forms.  The terminology is only shorthand for this equation; it is not a cohomological assertion.

There is a point which is essential in the calculus of integral forms.  The operator \(D_V\) is not, in general, a square-zero differential on the full algebra of integral forms.  Indeed, for an odd vector field the contraction \(i_V\) acts non-trivially on the delta symbols and their derivatives, and \(i_V^2\) need not vanish on arbitrary integral forms.  For example, already for the local field \(\partial/\partial\xi\), contraction sends \(\delta(d\xi)\) to a derivative delta symbol, and a second contraction gives a second derivative delta symbol.  Thus one must not justify localisation by the formal assertion \(D_V^2=0\).

The deformation argument uses only the following weaker identity.  If \(\omega\) is an ordinary one-form, then \(i_V\omega\) is a function, hence \(i_V(i_V\omega)=0\).  Moreover, since \(V\) is holomorphic, \(\ddbar\) and \(i_V\) super-commute on ordinary forms.  Therefore
\begin{align}
D_V(D_V\omega)
&=(\ddbar+i_V)(\ddbar\omega+i_V\omega)\notag\\
&=\ddbar^2\omega+\ddbar(i_V\omega)+i_V(\ddbar\omega)+i_V(i_V\omega)\notag\\
&=0.
\label{eq:DV-localising-square}
\end{align}
This identity is asserted only for one-forms \(\omega\).  It is not a replacement for a false statement that \(D_V^2\) vanishes on arbitrary integral forms.

\begin{definition}\label{def:localising}
Let \(U\subset S\) be open.  A smooth one-form \(\omega\) on \(U\) is \emph{\(V\)-localising} on a closed subset \(K\subset U_{\red}\) if \((i_V\omega)_{\red}\) is strictly positive on \(K\).  Near an isolated reduced zero, this positivity is required on the punctured reduced neighbourhood.
\end{definition}

\begin{lemma}\label{lem:localising-form}
Assume that, near a reduced zero, \(V\) has triangular form
\begin{equation}\label{eq:triangular-local-form-prelim}
V=\sum_{i=1}^n f_i(z,\xi)\frac{\partial}{\partial z_i}
 +\sum_{i=1}^n g_i(z)\frac{\partial}{\partial\xi_i},
\qquad \det J_g(0)\neq0.
\end{equation}
With respect to the standard Hermitian metric on the coordinate ball, put
\begin{equation}\label{eq:omega-local}
\omega=\sum_{i=1}^n\overline{g_i}\,dz_i+
       \sum_{i=1}^n\overline{g_i}\,d\xi_i.
\end{equation}
Then \(\omega\) is \(V\)-localising on a sufficiently small punctured neighbourhood of the zero, and
\begin{equation}\label{eq:positive-contraction}
(i_V\omega)_{\red}=\sum_{i=1}^n |g_i|^2.
\end{equation}
\end{lemma}

\begin{proof}
Contraction gives
\begin{equation}\label{eq:contraction-localising}
 i_V\omega=\sum_i\overline{g_i}f_i+
             \sum_i\overline{g_i}g_i.
\end{equation}
The first summand has zero reduction, since each \(f_i\) is odd.  The second summand reduces to \(\sum_i |g_i|^2\).  Since \(\det J_g(0)\ne0\), the holomorphic map \(g=(g_1,\ldots,g_n)\) has an isolated zero at the origin.  Therefore \(\sum_i |g_i|^2\) is strictly positive on a sufficiently small punctured reduced neighbourhood.
\end{proof}

\section{The global deformation argument}\label{sec:global}

The proof is a Dolbeault--Cartan deformation argument.  It uses the Stokes theorem for compactly supported integral forms on complex supermanifolds and the trace property of the Berezin integral, namely the vanishing of total odd derivatives \cite{Berezin,Manin,BernsteinLeites,Voronov,WittenNotes}.

Let \(\mathscr I^{\bullet}_S\) denote the sheaf of smooth integral forms.  Integration of a maximal-picture representative is defined locally by projection to the Berezinian coefficient, followed by ordinary integration on \(S_{\red}\).  The following lemma is the precise form of Stokes' theorem needed for the localisation operator.

\begin{lemma} \label{lem:exact-zero}
Let \(\alpha\in\Gamma_c(S,\mathscr I_S^{\bullet})\) be a compactly supported smooth integral form such that \(D_V\alpha\) has maximal picture and Berezinian degree.  Then
\[
\int_S D_V\alpha=0.
\]
\end{lemma}

\begin{proof}
By choosing a partition of unity on \(S_{\red}\), the assertion reduces to the corresponding statement in a coordinate chart.  Write local coordinates as
\[
(z_1,\ldots,z_n;\xi_1,\ldots,\xi_n)
\]
and write a compactly supported integral form as a finite sum of monomials
\[
F(z,\bar z,\xi,\bar\xi)\,
 dz_I\wedge d\bar z_J\,
 \delta^{(\nu)}(d\xi)\delta^{(\mu)}(d\bar\xi),
\]
with smooth compactly supported coefficients.  The contribution of \(\bar\partial\alpha\) to the Berezinian coefficient is a finite sum of ordinary \(\bar\partial\)-derivatives on the reduced manifold.  Hence its integral over \(S_{\red}\) is zero by the ordinary Stokes theorem.

It remains to examine the contraction term.  Locally
\[
V=\sum_i f_i\frac{\partial}{\partial z_i}+\sum_i g_i\frac{\partial}{\partial\xi_i}.
\]
The contraction \(i_{\partial/\partial z_i}\) lowers the ordinary holomorphic degree.  A monomial which does not recover the Berezinian degree after this operation has zero projection to the basic symbol.  In the remaining monomials the coefficient obtained after projection is a compactly supported function on \(S_{\red}\).  The odd contraction is the standard trace-zero operation on the exterior algebra of the odd coordinates; after projection to the Berezinian coefficient it is a finite sum of terms of the form
\[
[\xi_1\cdots\xi_n\bar\xi_1\cdots\bar\xi_n]
\frac{\partial H}{\partial\xi_i}
\]
or the corresponding expression with \(\bar\xi_i\).  This coefficient is zero, because Berezin integration is the trace on the exterior algebra and annihilates superderivatives.  In the integral-form notation this is the same statement as the distributional identity
\[
d\xi_i\delta^{(r)}(d\xi_i)=-r\delta^{(r-1)}(d\xi_i),
\]
which produces no boundary term in the odd variables.  Thus the integral of the \(i_V\alpha\)-contribution is zero.  Summing the \(\bar\partial\)- and \(i_V\)-contributions proves the lemma.
\end{proof}

\begin{lemma}[Deformation invariance]\label{lem:deformation-invariance}
Let \(\eta\) be a compactly supported \(D_V\)-closed integral form of maximal picture, and let \(\omega\) be a smooth one-form on a neighbourhood of \(\supp(\eta)\).  Assume that the product
\[
\eta\exp\left(-\frac{D_V\omega}{t}\right)
\]
is interpreted in the finite algebra of integral forms.  Then
\[
I(t)=\int_S\eta\exp\left(-\frac{D_V\omega}{t}\right),\qquad t>0,
\]
is independent of \(t\).
\end{lemma}

\begin{proof}
By \eqref{eq:DV-localising-square}, one has
$
D_V(D_V\omega)=0.
$
After multiplication by a fixed integral-form representative and projection to the Berezinian coefficient, the exponential is finite: the nilpotent part of \(i_V\omega\) is nilpotent, and multiplication by the integral-form part lowers only finitely many delta-derivative orders.  Therefore differentiation under the integral sign is legitimate.  The graded Leibniz rule gives
\begin{align*}
D_V\left(\omega\exp\left(-\frac{D_V\omega}{t}\right)\right)
&=(D_V\omega)\exp\left(-\frac{D_V\omega}{t}\right)  \\
&\quad -\omega\frac{1}{t}D_V(D_V\omega)\exp\left(-\frac{D_V\omega}{t}\right) \\
&=(D_V\omega)\exp\left(-\frac{D_V\omega}{t}\right).
\end{align*}
Since \(D_V\eta=0\), again by the graded Leibniz rule,
\begin{align*}
\frac{dI}{dt}
&=\frac{1}{t^2}
\int_S\eta(D_V\omega)\exp\left(-\frac{D_V\omega}{t}\right)\\
&=\frac{1}{t^2}
\int_S D_V\left(\eta\omega\exp\left(-\frac{D_V\omega}{t}\right)\right).
\end{align*}
The last integral is zero by Lemma~\ref{lem:exact-zero}.  Hence \(I(t)\) is constant.
\end{proof}

\begin{lemma} \label{lem:global-localising-form}
Let \(K\subset S_{\red}\) be compact and suppose that the reduced section \(s_V\in H^0(S_{\red},\cE)\) has no zero on \(K\).  Then there is a smooth one-form \(\omega\) on a neighbourhood of \(K\) such that
\[
(i_V\omega)_{\red}>0
\quad\text{on }K.
\]
\end{lemma}

\begin{proof}
Choose a Hermitian metric \(h\) on the vector bundle \(\cE\) over a neighbourhood of \(K\).  Since \(s_V\) has no zero on \(K\), the smooth section
\[
s_V^\flat:=h(s_V,\cdot)\in \Gamma(K,\cE^\vee)
\]
satisfies \(s_V^\flat(s_V)=\|s_V\|_h^2>0\) on \(K\).  The restriction of the odd cotangent sheaf to the reduced manifold maps onto \(\cE^\vee\).  After shrinking the neighbourhood and using a smooth partition of unity, choose a smooth one-form \(\omega\) whose restriction to the odd cotangent quotient is \(s_V^\flat\).  Then the reduced scalar part of \(i_V\omega\) is precisely \(s_V^\flat(s_V)=\|s_V\|_h^2\), and is therefore strictly positive on \(K\).
\end{proof}

\begin{theorem}[Vanishing]\label{thm:vanishing}
Let \(S\) be a compact complex supermanifold of dimension \(n|n\), and let \(V\) be an odd holomorphic vector field whose reduced zero locus is empty.  If \(\eta\) is a smooth \(D_V\)-closed integral form of maximal picture, then
\[
\int_S\eta=0.
\]
\end{theorem}

\begin{proof}
By Lemma~\ref{lem:global-localising-form}, applied to \(K=S_{\red}\), there is a smooth one-form \(\omega\) such that
\[
\rho:=(i_V\omega)_{\red}>0
\quad\text{on }S_{\red}.
\]
Compactness gives a constant \(c>0\) with \(\rho\geq c\).  By Lemma~\ref{lem:deformation-invariance}, the function
\[
I(t)=\int_S\eta\exp\left(-\frac{D_V\omega}{t}\right)
\]
is independent of \(t>0\).  Since
\[
\lim_{t\to\infty}\exp\left(-\frac{D_V\omega}{t}\right)=1
\]
in the finite integral-form algebra, one has
\[
I(t)=\lim_{T\to\infty}I(T)=\int_S\eta.
\]
On the other hand, writing
$
D_V\omega=\rho+R,
$
where \(R\) is a sum of nilpotent and positive form-degree terms, the exponential has the form
\[
\exp\left(-\frac{D_V\omega}{t}\right)
=e^{-\rho/t}\sum_{r=0}^N t^{-r}P_r,
\]
with fixed smooth integral forms \(P_r\).  The integer \(N\) is finite because the odd variables and the integral-form delta derivatives are finite in every local representative.  Since \(e^{-\rho/t}\leq e^{-c/t}\), every coefficient of the integral form of maximal picture is bounded by a finite sum of terms
\[
C_r t^{-r}e^{-c/t}.
\]
These tend to zero as \(t\to0^+\) and are dominated on compact sets.  Therefore
\[
\lim_{t\to0^+}I(t)=0.
\]
Since \(I(t)=\int_S\eta\), the desired vanishing follows.
\end{proof}

\begin{definition}\label{def:local-residue}
Let \(p\in Z(V)_{\red}\) be an isolated reduced zero, and let \(\eta\) be a smooth maximal-picture integral form defined near \(p\).  Choose a coordinate neighbourhood \(U\) of \(p\), a smooth cut-off \(\chi\) with compact support in \(U\) and equal to one on a smaller neighbourhood of \(p\), and a smooth one-form \(\omega_p\) on \(U\) such that
\[
(i_V\omega_p)_{\red}>0
\quad\text{on }U_{\red}\setminus\{p\}.
\]
When the following limit is evaluated by the finite integral-form expansion described above, set
\begin{equation}\label{eq:local-residue}
\Res_p(V,\eta)=
\lim_{t\to0^+}
\int_U \chi\eta\exp\left(-\frac{D_V\omega_p}{t}\right).
\end{equation}
Only the germ of \(\eta\) at \(p\) is meant.  If \(\eta\) is already compactly supported in \(U\) and the cut-off is identically one on its support, this reduces to the preceding compactly supported expression.
\end{definition}

\begin{lemma} \label{lem:residue-independent}
If \(D_V\eta=0\) in a neighbourhood of \(p\), then \(\Res_p(V,\eta)\) is independent of the cut-off \(\chi\), of the sufficiently small coordinate neighbourhood \(U\), and of the localising one-form \(\omega_p\).
\end{lemma}

\begin{proof}
First fix \(\omega_p\).  If \(\chi_0\) and \(\chi_1\) are two cut-offs equal to one near \(p\), then \(\chi_1-\chi_0\) is supported in an annulus on which \((i_V\omega_p)_{\red}\geq c>0\).  Hence the same estimate as in Theorem~\ref{thm:vanishing} gives
\[
\lim_{t\to0^+}
\int_U (\chi_1-\chi_0)\eta
\exp\left(-\frac{D_V\omega_p}{t}\right)=0.
\]
Thus the residue is independent of the cut-off and, consequently, of shrinking \(U\).
It remains to vary the localising form.  Let \(\omega_0\) and \(\omega_1\) be two choices.  After shrinking the neighbourhood and fixing a cut-off \(\chi\) supported there, the convex combination
\[
\omega_s=(1-s)\omega_0+s\omega_1,
\qquad 0\leq s\leq1,
\]
has reduced scalar part uniformly positive on the annulus where \(d\chi\neq0\) and on the boundary of the support of \(\chi\).  Differentiating the cut-off local integral gives
\[
\frac{\partial}{\partial s}
\int_U\chi\eta e^{-D_V\omega_s/t}
=-\frac{1}{t}
\int_U\chi\eta D_V(\dot\omega_s)e^{-D_V\omega_s/t}.
\]
Since \(D_V\eta=0\) near \(p\), the right hand side differs from
\[
-\frac{1}{t}
\int_U D_V\bigl(\chi\eta\dot\omega_s e^{-D_V\omega_s/t}\bigr)
\]
only by terms containing \(D_V\chi\).  The total \(D_V\)-term has zero integral by Lemma~\ref{lem:exact-zero}; the remaining \(D_V\chi\)-term is supported where the reduced scalar part of \(D_V\omega_s\) is bounded below by a positive constant.  Its contribution therefore tends to zero as \(t\to0^+\).  The limiting residue is independent of \(s\), hence of the localising form.
\end{proof}

\begin{theorem} \label{thm:residue-decomposition}
Let \(S\) be compact of dimension \(n|n\), and let \(V\) be an odd holomorphic vector field with finite reduced zero locus
$
Z(V)_{\red}=\{p_1,\ldots,p_N\}.
$
Then every \(D_V\)-closed smooth integral form of maximal picture \(\eta\) satisfies
\[
\int_S\eta=\sum_{j=1}^N\Res_{p_j}(V,\eta).
\]
\end{theorem}

\begin{proof}
Choose pairwise disjoint coordinate balls \(U_j\) around the points \(p_j\), and choose smaller balls \(U_j'\Subset U_j\).  On each \(U_j\) choose a localising one-form \(\omega_j\), smooth on \(U_j\), with \((i_V\omega_j)_{\red}>0\) on \(U_j\setminus\{p_j\}\).  On the compact set
\[
K=S_{\red}\setminus\bigcup_j U_j'
\]
the reduced section \(s_V\) has no zero; hence Lemma~\ref{lem:global-localising-form} gives a localising form \(\omega_0\) near \(K\).  Choose a smooth partition of unity
$
\chi_0+\chi_1+\cdots+\chi_N=1
$
on \(S_{\red}\), with \(\chi_j\) supported in \(U_j\) and equal to one near \(p_j\), and extend these functions to even smooth functions on \(S\).  Put
\[
\omega=\chi_0\omega_0+\sum_{j=1}^N\chi_j\omega_j.
\]
The reduced scalar part \((i_V\omega)_{\red}\) is strictly positive away from the points \(p_j\).  Indeed,
\[
(i_V\omega)_{\red}
=\chi_0(i_V\omega_0)_{\red}
+\sum_{j=1}^N\chi_j(i_V\omega_j)_{\red},
\]
and the right hand side is a non-negative sum with at least one positive summand at every point outside the chosen zeros.
By Lemma~\ref{lem:deformation-invariance},
\[
\int_S\eta=
\int_S\eta\exp\left(-\frac{D_V\omega}{t}\right)
\]
for every \(t>0\), using the limit as \(t\to\infty\) as in the proof of Theorem~\ref{thm:vanishing}.  Split the last integral by the same partition of unity.  On the support of \(\chi_0\), the reduced scalar part of \(D_V\omega\) is bounded below by a positive constant, and the same exponential estimate as in Theorem~\ref{thm:vanishing} shows that this contribution tends to zero as \(t\to0^+\).  For \(j\geq1\), the function \(\chi_j\) is equal to one near \(p_j\), and the part of its support outside a smaller neighbourhood of \(p_j\) again has positive reduced scalar part; hence that annular contribution tends to zero.  The remaining limit is exactly the cut-off definition of \(\Res_{p_j}(V,\eta)\).  Summing over \(j\) gives the formula.
\end{proof}

\section{Local calculation at a non-degenerate reduced zero}\label{sec:local}

The local computation separates two independent mechanisms.  The reduced variables produce an ordinary complex Gaussian integral and hence the holomorphic Jacobian.  The odd variables contribute through the finite algebra of integral forms and the lowering identities for delta symbols.  No exterior top degree in the even differentials \(d\xi_i\) is used.

Let the reduced zero be the origin.  The explicit formula is proved under the following normal-form hypothesis.

\begin{definition} \label{def:triangular-nondeg}
The odd holomorphic vector field \(V\) is said to be triangular and non-degenerate at the origin if there are local holomorphic supercoordinates \((z_1,\ldots,z_n;\xi_1,\ldots,\xi_n)\), centred at the reduced zero, such that
\begin{equation}\label{eq:triangular}
V=
\sum_{i=1}^n f_i(z,\xi)\frac{\partial}{\partial z_i}
+
\sum_{i=1}^n g_i(z)\frac{\partial}{\partial\xi_i},
\qquad
\det J_g(0)\neq0,
\end{equation}
where
\[
J_g=\left(\frac{\partial g_i}{\partial z_j}\right)_{i,j}.
\]
Thus the reduced zero equations are precisely
$\{
g_1(z)=\cdots=g_n(z)=0
\}$
on the reduced manifold.  This hypothesis is not a formal consequence of the existence of an isolated zero; it is the normal form under which the explicit formula below is proved.
\end{definition}

Let
\begin{equation}\label{eq:omega-split-expanded}
\omega=\omega_1+\omega_2,
\qquad
\omega_1=\sum_{i=1}^n\overline{g_i}\,dz_i,
\qquad
\omega_2=\sum_{i=1}^n\overline{g_i}\,d\xi_i.
\end{equation}
Then
\begin{equation}\label{eq:DVomega-split-expanded}
D_V\omega=A+B+C,
\end{equation}
where
\begin{align}
A&=i_V\omega,
\label{eq:A-expanded}\\
B&=\ddbar\omega_1
 =\sum_{i,j}\overline{\frac{\partial g_i}{\partial z_j}}\,d\bar z_j\wedge dz_i,
\label{eq:B-expanded}\\
C&=\ddbar\omega_2
 =\sum_{i,j}\overline{\frac{\partial g_i}{\partial z_j}}\,d\bar z_j\wedge d\xi_i.
\label{eq:C-expanded}
\end{align}
The scalar part is
\begin{equation}\label{eq:A-expanded-reduced}
A_{\red}=\sum_{i=1}^n |g_i|^2.
\end{equation}
We shall write
\begin{equation}\label{eq:Anil-expanded}
A_{\nil}:=A-A_{\red}=\sum_i\overline{g_i}f_i.
\end{equation}
Indeed,
\[
i_V\omega=
\sum_i\overline{g_i}f_i+
\sum_i\overline{g_i}g_i,
\]
and the first summand has zero reduction because the \(f_i\) are odd.  Notice, however, that \(A_{\nil}\) must not be ignored in the local expansion: although it is nilpotent in the odd variables, it is multiplied by powers of \(t^{-1}\) in the exponential and can contribute to the normalised Gaussian limit unless an additional vanishing condition is imposed.

\subsection{The ordinary differential part}

The holomorphic determinant enters through the following identities.

\begin{lemma}\label{lem:powers-B-expanded}
With \(B\) as in \eqref{eq:B-expanded}, one has
\begin{align}
B^{n-1}
&=(n-1)!\sum_{k=1}^n
\bigwedge_{i\neq k}
\left(d\overline{g_i}\wedge dz_i\right),
\label{eq:B-nminus-one-expanded}\\
B^n
&=n!\,d\overline{g_1}\wedge dz_1\wedge\cdots\wedge
 d\overline{g_n}\wedge dz_n
\label{eq:B-n-first-expanded}\\
&=(-1)^{n(n+1)/2}n!\det(\overline{J_g})
\,dz_1\wedge\cdots\wedge dz_n
\wedge d\bar z_1\wedge\cdots\wedge d\bar z_n .
\label{eq:B-n-expanded}
\end{align}
\end{lemma}

\begin{proof}
Since
$
B=\sum_i d\overline{g_i}\wedge dz_i,
$
the only non-zero terms in \(B^{n-1}\) are those in which all but one of the factors \(d\overline{g_i}\wedge dz_i\) occur.  This gives \eqref{eq:B-nminus-one-expanded}.  Similarly,
\[
B^n=n!\,d\overline{g_1}\wedge dz_1\wedge\cdots\wedge d\overline{g_n}\wedge dz_n.
\]
Moving each \(d\overline{g_i}\) past the holomorphic differentials gives the sign
$
(-1)^{1+2+\cdots+n}=(-1)^{n(n+1)/2}.
$
Finally,
\[
d\overline{g_i}=\sum_{j=1}^n
\overline{\frac{\partial g_i}{\partial z_j}}\,d\bar z_j,
\]
and hence
$
d\overline{g_1}\wedge\cdots\wedge d\overline{g_n}
=\det(\overline{J_g})\,d\bar z_1\wedge\cdots\wedge d\bar z_n.
$
Substitution proves \eqref{eq:B-n-expanded}.
\end{proof}

\begin{lemma}\label{lem:BC-expanded}
The mixed product \(B^{n-1}C\) is
\begin{equation}\label{eq:B-nminus-one-C-expanded}
B^{n-1}C
=(n-1)!\sum_{k=1}^n
\left(
\bigwedge_{i\neq k} d\overline{g_i}\wedge dz_i
\right)
\wedge d\overline{g_k}\wedge d\xi_k .
\end{equation}
Equivalently, after moving the holomorphic differentials to the front,
\begin{equation}\label{eq:B-nminus-one-C-reordered}
B^{n-1}C
=(-1)^{n(n+1)/2}(n-1)!\sum_{k=1}^n
 dz_1\wedge\cdots\wedge d\xi_k\wedge\cdots\wedge dz_n
\wedge d\overline{g_1}\wedge\cdots\wedge d\overline{g_n}.
\end{equation}
\end{lemma}

\begin{proof}
In the product \(B^{n-1}C\), the factor \(C\) is
\[
C=\sum_k d\overline{g_k}\wedge d\xi_k.
\]
For a fixed \(k\), the contribution which can contain the full anti-holomorphic determinant is obtained by taking from \(B^{n-1}\) all factors \(d\overline{g_i}\wedge dz_i\) with \(i\neq k\), and from \(C\) the factor \(d\overline{g_k}\wedge d\xi_k\).  If any index is repeated, two anti-holomorphic factors become linearly dependent and the contribution to the determinant term vanishes.  This gives \eqref{eq:B-nminus-one-C-expanded}.  The reordering is the same as in Lemma~\ref{lem:powers-B-expanded}, except that the place of \(dz_k\) is occupied by \(d\xi_k\), giving \eqref{eq:B-nminus-one-C-reordered}.
\end{proof}

\begin{remark}\label{rem:no-C-square}
The calculation does not use any assertion of the form \(C^2=0\).  In the usual calculus of supermanifolds, the differentials \(d\xi_i\) are even, and there is no exterior top degree in these variables.  The finiteness of the local expression comes instead from the integral-form identities
\[
d\xi_i\delta^{(r)}(d\xi_i)=-r\delta^{(r-1)}(d\xi_i),
\]
which lower the order of the delta derivative.  Thus, for a fixed integral-form representative, only finitely many powers of \(C\) can contribute to the coefficient of the basic Berezinian symbol.
\end{remark}

\subsection{The Gaussian normalisation}

\begin{lemma} \label{lem:bosonic-gaussian}
For every smooth compactly supported function \(\varphi\) on a sufficiently small reduced coordinate ball,
\begin{equation}\label{eq:gaussian-expanded}
\lim_{t\to0^+}
\int
\varphi(z,\bar z)e^{-\sum_i|g_i(z)|^2/t}
\frac{(-1)^nB^n}{t^n n!}
=
\left(\frac{2\pi}{i}\right)^n
\frac{\varphi(0)}{\det J_g(0)}.
\end{equation}
\end{lemma}

\begin{proof}
By the holomorphic inverse function theorem, \(u_i=g_i(z)\) is a holomorphic coordinate system near the origin.  This is the same analytic mechanism underlying the Grothendieck--Bott residue denominator \cite{Bott,CarrellLieberman,Liu}.  In these coordinates,
\[
\sum_i|g_i(z)|^2=\sum_i|u_i|^2.
\]
Moreover,
\[
du_1\wedge\cdots\wedge du_n=\det J_g(z)\,dz_1\wedge\cdots\wedge dz_n,
\]
and
\[
d\bar u_1\wedge\cdots\wedge d\bar u_n=\det(\overline{J_g(z)})\,d\bar z_1\wedge\cdots\wedge d\bar z_n.
\]
Using Lemma~\ref{lem:powers-B-expanded}, together with the factor \((-1)^n\) displayed in \eqref{eq:gaussian-expanded}, the anti-holomorphic determinant appearing in \(B^n\) cancels the anti-holomorphic Jacobian of the change of variables, with the chosen orientation convention for \(du_1\wedge\cdots\wedge du_n\wedge d\bar u_1\wedge\cdots\wedge d\bar u_n\).  The holomorphic Jacobian remains in the denominator.  Hence the left-hand side of \eqref{eq:gaussian-expanded} is equal to
\[
\lim_{t\to0^+}
\int
\frac{\varphi(z(u),\bar z(\bar u))}{\det J_g(z(u))}
 e^{-|u|^2/t}
\frac{du_1\wedge\cdots\wedge du_n\wedge d\bar u_1\wedge\cdots\wedge d\bar u_n}{t^n}.
\]
The standard complex Gaussian identity
\[
\int_{\C^n}
 e^{-|u|^2/t}
\frac{du_1\wedge\cdots\wedge du_n\wedge d\bar u_1\wedge\cdots\wedge d\bar u_n}{t^n}
=
\left(\frac{2\pi}{i}\right)^n
\]
then gives \eqref{eq:gaussian-expanded}.
\end{proof}

\section{Integral-form coefficient extraction}\label{sec:coefficient-extraction}

The coefficient calculation is purely algebraic once the reduced Gaussian factor has been separated.  The integration map used throughout this paper is the projection onto the ordinary top degree and onto the basic integral-form symbol.  The exponential factor in the localisation argument has a precise role: its ordinary differential part supplies the missing reduced differentials, while its integral-form part lowers derivative delta symbols until the basic symbol is reached.

Set
\begin{equation}\label{eq:Delta-basic-v13}
\Delta=
\delta(d\xi_1)\cdots\delta(d\xi_n)
\delta(d\bar\xi_1)\cdots\delta(d\bar\xi_n)
\end{equation}
and
\begin{equation}\label{eq:D-basic-again}
\mathcal D=dz_1\wedge\cdots\wedge dz_n\wedge d\bar z_1\wedge\cdots\wedge d\bar z_n\,\Delta .
\end{equation}
Let
\[
\Pi_{\mathcal D}
\]
denote the projection of a local integral form onto the coefficient of \(\mathcal D\).  More explicitly, if \(\Theta\) is a finite sum of terms involving ordinary differentials and derivatives of the delta symbols, then \(\Pi_{\mathcal D}(\Theta)\) is the unique smooth superfunction \(H\), if it exists, such that the component of \(\Theta\) proportional to \(\mathcal D\) is \(H\mathcal D\).  The local integral is obtained by applying, after this projection, the Berezin coefficient extraction
\[
H\longmapsto
[\xi_1\cdots\xi_n\bar\xi_1\cdots\bar\xi_n]H
\]
and then integrating the resulting ordinary top form on the reduced manifold.

Equivalently, \(\Pi_{\mathcal D}\) is obtained by the following explicit procedure.  A monomial is first reordered to the normal form
\[
F\,dz_1\cdots dz_n\,d\bar z_1\cdots d\bar z_n
\prod_i (d\xi_i)^{a_i}\delta^{(\nu_i)}(d\xi_i)
\prod_i (d\bar\xi_i)^{b_i}\delta^{(\mu_i)}(d\bar\xi_i),
\]
with the Koszul sign produced by the reordering included in \(F\).  The products in the even variables are then reduced by
\[
(d\xi_i)^a\delta^{(r)}(d\xi_i)
=\begin{cases}
(-1)^a\dfrac{r!}{(r-a)!}\delta^{(r-a)}(d\xi_i),& a\le r,\\[0.4em]
0,& a>r,
\end{cases}
\]
and by the analogous formula for \(d\bar\xi_i\).  The monomial contributes to \(\Pi_{\mathcal D}\) precisely when all final delta orders are zero and no ordinary differential is missing or repeated.

A general local representative of maximal picture can be written as a finite sum
\begin{equation}\label{eq:general-integral-representative}
\eta=
\sum_{I,J,\nu,\mu}
\Phi_{I,J,\nu,\mu}(z,\bar z,\xi,\bar\xi)
\,dz_I\wedge d\bar z_J\,
\delta^{(\nu)}(d\xi)\delta^{(\mu)}(d\bar\xi),
\end{equation}
where \(I,J\subset\{1,\ldots,n\}\), the functions \(\Phi_{I,J,\nu,\mu}\) are smooth superfunctions, and
\[
\delta^{(\nu)}(d\xi)=
\delta^{(\nu_1)}(d\xi_1)\cdots\delta^{(\nu_n)}(d\xi_n),
\qquad
\delta^{(\mu)}(d\bar\xi)=
\delta^{(\mu_1)}(d\bar\xi_1)\cdots\delta^{(\mu_n)}(d\bar\xi_n).
\]
The expansion is finite.  No exterior top degree in the variables \(d\xi_i\) is being assumed.

\begin{lemma} \label{lem:lowering-rule-expanded}
For every \(r\geq0\),
\begin{equation}\label{eq:lowering-rule-expanded}
d\xi_i\,\delta^{(r)}(d\xi_i)=
-r\,\delta^{(r-1)}(d\xi_i),
\end{equation}
where the right hand side is interpreted as zero for \(r=0\).  More generally,
\begin{equation}\label{eq:iterated-lowering-rule-expanded}
(d\xi_i)^a\delta^{(r)}(d\xi_i)
=\begin{cases}
(-1)^a\dfrac{r!}{(r-a)!}\delta^{(r-a)}(d\xi_i),& a\le r,\\[0.4em]
0,& a>r.
\end{cases}
\end{equation}
The same formula holds with \(\xi_i\) replaced by \(\bar\xi_i\).  Consequently, multiplication by
\[
C=\sum_{i,j}\overline{\frac{\partial g_i}{\partial z_j}}\,d\bar z_j\wedge d\xi_i
\]
lowers exactly one holomorphic delta order \(\nu_i\) by one and adds one anti-holomorphic ordinary differential \(d\bar z_j\), unless \(\nu_i=0\), in which case that term vanishes.
\end{lemma}

\begin{proof}
The first identity is the standard distributional identity for integral forms \cite{BernsteinLeites,Voronov,WittenNotes}.  It is obtained by differentiating the relation \(d\xi_i\delta(d\xi_i)=0\) with respect to the even variable \(d\xi_i\).  Iterating this identity gives
\[
(d\xi_i)^a\delta^{(r)}(d\xi_i)
=(-r)(-(r-1))\cdots (-(r-a+1))\delta^{(r-a)}(d\xi_i)
\]
for \(a\le r\), which is \eqref{eq:iterated-lowering-rule-expanded}.  If \(a>r\), the iteration reaches \(d\xi_i\delta(d\xi_i)=0\).  Since the variables \(d\xi_i\) are even, this is the mechanism which replaces exterior nilpotence.  Applying \eqref{eq:lowering-rule-expanded} to each summand of \(C\) gives the stated effect.  The finiteness follows from the fact that, in \eqref{eq:general-integral-representative}, only finitely many derivative orders occur.
\end{proof}

Write
$
A=A_{\red}+A_{\nil},
$ and $
A_{\red}=\sum_i|g_i|^2.
$
For coefficient extraction one uses the ordered exponential expansion
\begin{equation}\label{eq:finite-expansion-full}
\exp\left(-\frac{A+B+C}{t}\right)
=e^{-A_{\red}/t}
\sum_{\ell\geq0}\frac{(-1)^\ell}{\ell!t^\ell}
(A_{\nil}+B+C)^\ell .
\end{equation}
The power in \eqref{eq:finite-expansion-full} is the product in the integral-form algebra; it is expanded as a sum of ordered words in the three letters \(A_{\nil}\), \(B\), and \(C\).  Only after a word has been multiplied by the fixed representative \(\eta\) do we normal-order it according to Convention~\ref{conv:normal-ordering}.  Thus no multinomial rearrangement is being assumed.  The projected sum is finite: the ordinary part must have bidegree \((n,n)\), and the integral-form part must be reduced to the undifferentiated symbol \(\Delta\).  Powers of \(C\) are not killed by exterior nilpotence; they stop contributing only when there are no remaining derivative delta symbols to lower.

\begin{definition} \label{def:wick-functional}
Let \(F(u,\bar u,\xi,\bar\xi)\) be a smooth superfunction near the origin of \(\C^n\).  Write its Taylor expansion in the reduced variables as
\[
F(u,\bar u,\xi,\bar\xi)
\sim
\sum_{\alpha,\beta}
F_{\alpha\beta}(\xi,\bar\xi)u^\alpha\bar u^\beta .
\]
For \(q\geq0\) define
\begin{equation}\label{eq:wick-functional-definition}
\mathfrak W_q(F)=
\sum_{|\alpha|=q}\alpha!\,F_{\alpha\alpha}(\xi,\bar\xi),
\end{equation}
and set \(\mathfrak W_q(F)=0\) for \(q<0\).  Equivalently,
\begin{equation}\label{eq:wick-functional-derivatives}
\mathfrak W_q(F)=
\sum_{|\alpha|=q}
\frac{1}{\alpha!}
\left.
\frac{\partial^{2|\alpha|}F}
{\partial u^\alpha\partial\bar u^\alpha}
\right|_{u=0} .
\end{equation}
Thus \(\mathfrak W_0(F)=F(0,0,\xi,\bar\xi)\), but for \(q>0\) the functional is not evaluation at the origin.  It records the diagonal Gaussian moments.
\end{definition}

\begin{lemma} \label{lem:wick-limit}
With the orientation convention fixed above, if \(F\) is represented by a smooth compactly supported superfunction on a sufficiently small reduced coordinate ball, then for every \(q\geq0\) one has
\begin{equation}\label{eq:wick-asymptotic}
\lim_{t\to0^+}
\int_{\C^n}
 t^{-n-q}F(u,\bar u,\xi,\bar\xi)e^{-|u|^2/t}
 du_1\wedge\cdots\wedge du_n\wedge d\bar u_1\wedge\cdots\wedge d\bar u_n
=
\left(\frac{2\pi}{i}\right)^n\mathfrak W_q(F).
\end{equation}
\end{lemma}

\begin{proof}
It is enough to check monomials.  If \(\alpha\neq\beta\), the integral of \(u^\alpha\bar u^\beta e^{-|u|^2/t}\) vanishes by the angular integrations.  If \(\alpha=\beta\), the product of the one-dimensional complex Gaussian moment identities gives
\[
\int_{\C^n}u^\alpha\bar u^\alpha e^{-|u|^2/t}
 du_1\wedge\cdots\wedge du_n\wedge d\bar u_1\wedge\cdots\wedge d\bar u_n
=
\left(\frac{2\pi}{i}\right)^n\alpha!\,t^{n+|\alpha|} .
\]
After multiplication by \(t^{-n-q}\), the limit is zero unless \(|\alpha|=q\), and is then the coefficient displayed in \eqref{eq:wick-functional-definition}.  The derivative expression \eqref{eq:wick-functional-derivatives} is the same formula written in Taylor coefficients.
\end{proof}

\begin{definition} \label{def:local-numerator-expanded}
Let \(\eta\) be a fixed local integral-form representative.  Expand, in the ordered sense of \eqref{eq:finite-expansion-full},
\begin{equation}\label{eq:numerator-expanded-expression}
\Pi_{\mathcal D}\left(
\eta\exp\left(-\frac{A+B+C}{t}\right)\right)
=
e^{-A_{\red}/t}\sum_{r=0}^N t^{-r}H_r(z,\bar z,\xi,\bar\xi)\mathcal D .
\end{equation}
Let \(u=g(z)\) be the holomorphic change of reduced variables near the non-degenerate zero, and let \(z=z(u)\) be its local inverse.  After changing variables in the ordinary reduced integral, set
\begin{equation}\label{eq:hat-H-definition}
\widehat H_r(u,\bar u,\xi,\bar\xi)=
H_r(z(u),\bar z(\bar u),\xi,\bar\xi)
\frac{\det J_g(0)}
{\det J_g(z(u))\det\overline{J_g(z(u))}} .
\end{equation}
The local numerator is
\begin{equation}\label{eq:numerator-definition-explicit}
\mathcal N_0(V,\eta)=
\left[
\xi_1\cdots\xi_n\bar\xi_1\cdots\bar\xi_n
\right]
\sum_{r=0}^N
\mathfrak W_{r-n}(\widehat H_r).
\end{equation}
Terms with \(r<n\) vanish because \(\mathfrak W_{r-n}=0\).  Terms with \(r>n\) are evaluated by Gaussian moments; they are not obtained by evaluating \(H_r\) at the origin.  In particular, the nilpotent scalar term \(A_{\nil}\) and the Taylor expansion of the Jacobian factor in \eqref{eq:hat-H-definition} are retained whenever their Wick degree is correct.
\end{definition}

\begin{definition}[\(A_{\nil}\)-silent representatives]\label{def:Anil-silent}
A local representative \(\eta\) is called \emph{\(A_{\nil}\)-silent} if, in the ordered expansion \eqref{eq:finite-expansion-full}, every word containing at least one occurrence of \(A_{\nil}\) has zero contribution to the finite Wick sum defining \(\mathcal N_0(V,\eta)\) after applying \(\Pi_{\mathcal D}\).  This is automatic, for example, if the chosen triangular model has \(f_i=0\) for all \(i\), hence \(A_{\nil}=0\).  Without this condition the terms containing \(A_{\nil}\) are part of the numerator \(\mathcal N_0(V,\eta)\) and must not be suppressed.
\end{definition}

The two elementary contributions below are the building blocks of the coordinate formulae.  The first corresponds to a representative with no ordinary differentials and no derivative delta symbols; the ordinary top degree is supplied by \(B^n\).  It gives the familiar simple numerator only under the \(A_{\nil}\)-silent hypothesis.

\begin{proposition} \label{prop:basic-representative-expanded}
Consider 
\begin{equation}\label{eq:basic-zero-representative-v13}
\eta=\Phi(z,\bar z,\xi,\bar\xi)\,\Delta .
\end{equation}
Assume that this representative is \(A_{\nil}\)-silent in the sense of Definition~\ref{def:Anil-silent}.  Then
\begin{equation}\label{eq:basic-numerator-expanded}
\mathcal N_0(V,\eta)=
[\xi_1\cdots\xi_n\bar\xi_1\cdots\bar\xi_n]\,
\Phi(0,0,\xi,\bar\xi).
\end{equation}
\end{proposition}

\begin{proof}
We compute the projection onto \(\mathcal D\).  Since the representative \eqref{eq:basic-zero-representative-v13} has no ordinary differentials, the only way to obtain the ordinary degree \((n,n)\) without using the lowering operator \(C\) is to take the factor \(B^n/n!\) from the exponential.  By Lemma~\ref{lem:powers-B-expanded},
\[
\frac{B^n}{n!}
=(-1)^{n(n+1)/2}\det(\overline{J_g})
\,dz_1\wedge\cdots\wedge dz_n
\wedge d\bar z_1\wedge\cdots\wedge d\bar z_n .
\]
Moreover, every factor \(C\) contains some \(d\xi_i\).  Since the integral-form part of \eqref{eq:basic-zero-representative-v13} contains \(\delta(d\xi_i)\) and no derivative delta symbols, the identity \(d\xi_i\delta(d\xi_i)=0\) shows that all terms containing \(C\) vanish under \(\Pi_{\mathcal D}\).  Thus the only possible contribution to \(\Pi_{\mathcal D}\) is
\begin{equation}\label{eq:basic-proj-v13}
\Pi_{\mathcal D}\left(
\Phi\Delta\,e^{-A_{\red}/t}\frac{(-1)^n}{t^n}\frac{B^n}{n!}
\right).
\end{equation}
By the \(A_{\nil}\)-silent hypothesis, all ordered words in the exponential which contain \(A_{\nil}\) have zero contribution to the finite Wick sum after projection.  Hence the surviving contribution is obtained from the word \(B^n\), by evaluating the smooth coefficient at \(z=0\) and extracting the top Berezin coefficient.  This gives \eqref{eq:basic-numerator-expanded}.  The determinant and the Gaussian normalisation are accounted for in Theorem~\ref{thm:general-local-residue}; they are not part of the numerator.
\end{proof}

\begin{remark} \label{rem:Anil-corrections}
The hypothesis in Proposition~\ref{prop:basic-representative-expanded} is not cosmetic.  Since
\[
A_{\nil}=\sum_i\overline{g_i}f_i,
\]
a word containing \(a\) copies of \(A_{\nil}\) comes with a factor \(t^{-a}\) in the exponential.  After the change of variables \(u=g(z)\), Gaussian moments of the form
\[
t^{-a}\int u^\alpha\bar u^\beta e^{-|u|^2/t}\,du\,d\bar u
\]
can have a finite non-zero limit when \(|\alpha|=|\beta|=a\).  Therefore one cannot discard \(A_{\nil}\) merely because it vanishes on the reduced zero.  The general numerator \(\mathcal N_0(V,\eta)\) was defined through the Wick functionals precisely so that these finite corrections, when present, are included.
\end{remark}

The next statement records the first non-trivial lowering phenomenon.  The indices are aligned in the unique configuration which contributes to the determinant term: one holomorphic differential is already present, and one odd delta derivative is lowered by the mixed part of \(D_V\omega\).

\begin{proposition}[One lowering step]\label{prop:one-lowering-step}
Fix \(k\in\{1,\ldots,n\}\).  Suppose that a summand of \(\eta\) is
\begin{equation}\label{eq:one-lowering-summand}
\eta_k=
\Psi_k(z,\bar z,\xi,\bar\xi)\,dz_k\,
\delta'(d\xi_k)\prod_{i\neq k}\delta(d\xi_i)
\prod_i\delta(d\bar\xi_i).
\end{equation}
Then the only part of \(B^{n-1}C\) which can reduce \(\eta_k\) to the basic symbol \(\mathcal D\) is
\begin{equation}\label{eq:one-lowering-relevant-factor}
(n-1)!\left(\bigwedge_{i\neq k}d\overline{g_i}\wedge dz_i\right)
\wedge d\overline{g_k}\wedge d\xi_k .
\end{equation}
After projection to \(\mathcal D\), its contribution is
\begin{equation}\label{eq:one-lowering-projection}
\Pi_{\mathcal D}\!\bigl(\eta_k\,B^{n-1}C\bigr)
=
\varepsilon_k (n-1)!\det(\overline{J_g})\,
\Psi_k\,
\mathcal D,
\end{equation}
where \(\varepsilon_k\in\{\pm1\}\) is the sign determined by the fixed ordering in \eqref{eq:D-basic-again}.  Consequently the ordered word displayed above contributes to \(\mathcal N_0(V,\eta)\) through the Wick degree zero term, namely by the Berezin coefficient of \(\varepsilon_k\Psi_k(0,0,\xi,\bar\xi)\).  Higher words containing \(A_{\nil}\) are handled separately by Definition~\ref{def:local-numerator-expanded}.
\end{proposition}

\begin{proof}
We first determine which terms can survive.  The representative \(\eta_k\) already contains the holomorphic ordinary differential \(dz_k\).  To obtain the full holomorphic ordinary degree, the factor \(B^{n-1}\) must therefore supply precisely the differentials \(dz_i\) with \(i\neq k\); any occurrence of \(dz_k\) from \(B\) would make the wedge product zero.  Similarly, the only way to lower the integral-form part from \(\delta'(d\xi_k)\) to \(\delta(d\xi_k)\) is to use the factor \(d\xi_k\) from \(C\).  If \(d\xi_i\) with \(i\neq k\) occurs, it multiplies \(\delta(d\xi_i)\) and hence vanishes by \(d\xi_i\delta(d\xi_i)=0\).  Thus the only possible contribution is exactly \eqref{eq:one-lowering-relevant-factor}.

By the lowering rule,
\[
d\xi_k\delta'(d\xi_k)=-\delta(d\xi_k).
\]
Moving the factor \(d\xi_k\) to the position of \(\delta'(d\xi_k)\), and then moving the ordinary differentials into the fixed order of \(\mathcal D\), produces a sign \(\varepsilon_k\).  The value of this sign depends only on the ordering convention in \eqref{eq:D-basic-again}; once that ordering has been fixed, it is independent of \(\Psi_k\), of \(V\), and of the coordinate values.

It remains to compute the ordinary determinant.  The surviving ordinary factor is
\[
dz_k\wedge\left(\bigwedge_{i\neq k}d\overline{g_i}\wedge dz_i\right)
\wedge d\overline{g_k}.
\]
After reordering the holomorphic differentials as \(dz_1\wedge\cdots\wedge dz_n\) and the anti-holomorphic differentials as \(d\overline{g_1}\wedge\cdots\wedge d\overline{g_n}\), this becomes, up to the same fixed sign,
\[
\det(\overline{J_g})
\,dz_1\wedge\cdots\wedge dz_n
\wedge d\bar z_1\wedge\cdots\wedge d\bar z_n.
\]
Together with the lowered delta symbols this is precisely \(\mathcal D\), proving \eqref{eq:one-lowering-projection}.  The statement about \(\mathcal N_0(V,\eta)\) follows from Definition~\ref{def:local-numerator-expanded} after applying the same Gaussian limit as in Lemma~\ref{lem:bosonic-gaussian}.
\end{proof}

\begin{remark}\label{rem:non-aligned-lowering}
The alignment in Proposition~\ref{prop:one-lowering-step} is essential.  For example, a term containing \(dz_\ell\delta'(d\xi_k)\) with \(\ell\neq k\) cannot contribute to the determinant term through a single lowering step: either \(B^{n-1}\) repeats a holomorphic differential, or the anti-holomorphic determinant has two identical rows and one missing row.  Such terms may only contribute in higher lowering patterns, and only when the ordinary and integral-form degrees match the projection to \(\mathcal D\).
\end{remark}

\section{The local residue formula}\label{sec:local-residue-formula}

\begin{theorem}[Non-degenerate local residue]\label{thm:general-local-residue}
Under the triangular non-degeneracy hypothesis \eqref{eq:triangular}, every local integral-form representative of a \(D_V\)-closed integral form of maximal picture satisfies
\begin{equation}\label{eq:general-local-residue-final}
\Res_0(V,\eta)=
\left(\frac{2\pi}{i}\right)^n
\frac{\mathcal N_0(V,\eta)}{\det J_g(0)}.
\end{equation}
In particular, for an \(A_{\nil}\)-silent basic zero-form representative \(\eta=\Phi\Delta\),
\begin{equation}\label{eq:basic-residue-final-expanded}
\Res_0(V,\eta)=
\left(\frac{2\pi}{i}\right)^n
\frac{[\xi_1\cdots\xi_n\bar\xi_1\cdots\bar\xi_n]\Phi(0,0,\xi,\bar\xi)}{\det J_g(0)}.
\end{equation}
\end{theorem}

\begin{proof}
Choose the localising form \(\omega\) of \eqref{eq:omega-split-expanded}.  By Lemma~\ref{lem:localising-form}, the reduced scalar part of \(D_V\omega\) is \(\sum_i|g_i|^2\), which is positive on a punctured neighbourhood of the origin.  The local residue is therefore computed from
\[
\int \eta\exp\left(-\frac{A+B+C}{t}\right),
\]
with \(A,B,C\) as in \eqref{eq:A-expanded}--\eqref{eq:C-expanded}.
Split \(A=A_{\red}+A_{\nil}\), where
$
A_{\red}=\sum_i|g_i|^2.
$
The nilpotent part is finite.  Hence, after multiplication by the fixed integral-form representative of \(\eta\), the exponential is expanded as
\[
 e^{-A_{\red}/t}
 \sum_{\ell\ge0}\frac{(-1)^\ell}{\ell!t^\ell}(A_{\nil}+B+C)^\ell,
\]
where the last power is the product in the integral-form algebra.  It is a sum over ordered words in \(A_{\nil},B,C\); no multinomial commutation is used.  Only those words whose ordinary differential part has top degree and whose integral-form part normal-orders to \(\mathcal D\) survive under \(\Pi_{\mathcal D}\).  This is exactly the ordered coefficient-extraction rule of Definition~\ref{def:local-numerator-expanded}.

After projection one obtains \eqref{eq:numerator-expanded-expression}.  Changing reduced variables by \(u=g(z)\), the ordinary volume element contributes the factor
\[
\frac{1}{\det J_g(z(u))\det\overline{J_g(z(u))}}.
\]
We separate the constant holomorphic denominator by multiplying the remaining coefficient by \(\det J_g(0)\), as in \eqref{eq:hat-H-definition}.  Thus every term in the projected integral is a finite sum of expressions of the form
\[
\frac{1}{\det J_g(0)}
\int t^{-r}\widehat H_r(u,\bar u,\xi,\bar\xi)
 e^{-|u|^2/t}
 du_1\wedge\cdots\wedge du_n\wedge d\bar u_1\wedge\cdots\wedge d\bar u_n .
\]
By the Wick limit, Lemma~\ref{lem:wick-limit}, the finite part of this integral is
\[
\left(\frac{2\pi}{i}\right)^n
\frac{1}{\det J_g(0)}
\mathfrak W_{r-n}(\widehat H_r).
\]
Summing over the finite list of projected ordered words and then extracting the top Berezin coefficient gives precisely \(\mathcal N_0(V,\eta)\).  This proves
\[
\Res_0(V,\eta)=
\left(\frac{2\pi}{i}\right)^n
\frac{\mathcal N_0(V,\eta)}{\det J_g(0)}.
\]
If \(\eta=\Phi\Delta\), Proposition~\ref{prop:basic-representative-expanded} gives the stated numerator, and \eqref{eq:basic-residue-final-expanded} follows.
\end{proof}

\begin{corollary} \label{cor:global-explicit-residue}
Let \(S\) be a compact complex supermanifold of dimension \(n|n\), and let \(V\) be an odd holomorphic vector field whose reduced zero locus is finite.  Assume that every reduced zero is triangular non-degenerate.  Then every \(D_V\)-closed smooth integral form of maximal picture \(\eta\) satisfies
\begin{equation}\label{eq:global-explicit-residue}
\int_S\eta=
\left(\frac{2\pi}{i}\right)^n
\sum_{p\in Z(V)_{\red}}
\frac{\mathcal N_p(V,\eta)}{\det J_g(p)}.
\end{equation}
In particular, if near each zero \(p\) the contributing local component is the basic zero-form component \(\Phi_p\Delta_p\), then the local summand is
\begin{equation}\label{eq:global-basic-explicit-residue}
\left(\frac{2\pi}{i}\right)^n
\frac{[\xi_1\cdots\xi_n\bar\xi_1\cdots\bar\xi_n]\Phi_p(p,\xi,\bar\xi)}{\det J_g(p)}.
\end{equation}
\end{corollary}

\begin{proof}
The residue decomposition theorem gives
\(\int_S\eta=\sum_{p\in Z(V)_{\red}}\Res_p(V,\eta)\).  At each zero, the triangular non-degeneracy hypothesis permits the local calculation of Theorem~\ref{thm:general-local-residue}.  Substituting the local formula into the residue decomposition gives \eqref{eq:global-explicit-residue}.  The expression \eqref{eq:global-basic-explicit-residue} is the special case of \eqref{eq:basic-residue-final-expanded} for a basic integral-form representative.
\end{proof}

\section{Explicit local models}\label{sec:explicit-local-models}

The following local models make explicit the projection mechanism used in Sections~\ref{sec:local} and~\ref{sec:coefficient-extraction}.  The representatives below are maximal-picture integral forms; whenever their ordinary degree is not top, the missing ordinary differentials are supplied by the exponential of \(D_V\omega\).

\subsection{Dimension \texorpdfstring{\(1|1\)}{1|1}}
Let
\[
V=f(z,\xi)\frac{\partial}{\partial z}+g(z)\frac{\partial}{\partial\xi},
\qquad g(0)=0,
\qquad g'(0)\neq0,
\]
and set
\[
\omega=\bar g\,dz+\bar g\,d\xi.
\]
Then
\[
D_V\omega=A+B+C,
\qquad
A=|g|^2+A_{\nil},
\qquad
A_{\nil}=\bar g f,
\qquad
B=\overline{g'}\,d\bar z\wedge dz,
\qquad
C=\overline{g'}\,d\bar z\wedge d\xi .
\]
Let \(\Delta=\delta(d\xi)\delta(d\bar\xi)\).  For the zero-form representative
 $
\eta_0=\Phi(z,\bar z,\xi,\bar\xi)\Delta,
$
the projection to the basic symbol coming from the word \(B\) is obtained from the factor \(-B/t\) in the exponential.  Since \(d\xi\delta(d\xi)=0\), every word containing \(C\) and no derivative delta symbol vanishes after projection.  If the representative is \(A_{\nil}\)-silent, no word containing \(A_{\nil}\) contributes to the finite Wick sum.  Thus, under this explicit hypothesis,
\[
\Pi_{\mathcal D}\left(\eta_0\exp\left(-\frac{D_V\omega}{t}\right)\right)
=
\Phi e^{-|g|^2/t}\left(-\frac{B}{t}\right)
\Delta+o(1)
\]
for the purpose of the local limit.  The change of variables \(u=g(z)\) gives
\[
\Res_0(V,\eta_0)=
\frac{2\pi}{i}\,
\frac{[\xi\bar\xi]\Phi(0,0,\xi,\bar\xi)}{g'(0)}.
\]

Now consider the representative
\[
\eta_1=\Psi(z,\bar z,\xi,\bar\xi)\,dz\,\delta'(d\xi)\delta(d\bar\xi).
\]
Here the holomorphic ordinary differential is already present, and the missing anti-holomorphic differential is supplied by the word \(C\).  Normal ordering gives the following calculation.  Starting from the ordered product of the representative with the relevant factor of the exponential, one obtains
\[
\Psi\,dz\,\delta'(d\xi)\delta(d\bar\xi)\,
\overline{g'}\,d\bar z\wedge d\xi .
\]
Move the ordinary differentials to the left, with the Koszul sign prescribed by Convention~\ref{conv:normal-ordering}, and move the even variable \(d\xi\) next to \(\delta'(d\xi)\).  Then
\[
d\xi\delta'(d\xi)=-\delta(d\xi),
\]
so the normal-ordered result is
\[
\epsilon_{1|1}\,\overline{g'}\,\Psi\,
 dz\wedge d\bar z\,\delta(d\xi)\delta(d\bar\xi)
=
\epsilon_{1|1}\,\overline{g'}\,\Psi\,\mathcal D,
\]
where \(\epsilon_{1|1}\in\{\pm1\}\) is fixed once and for all by the chosen order of the basic symbol.  Thus the \(C/t\)-term contributes a finite numerator obtained from the Berezin coefficient of \(\epsilon_{1|1}\Psi\), together with any additional \(A_{\nil}\)-corrections and Jacobian Taylor corrections allowed by the Wick definition of Definition~\ref{def:local-numerator-expanded}.  This is the simplest instance of the lowering rule: \(C\) lowers a derivative delta symbol; it is not an exterior-nilpotent term.

\subsection{Dimension \texorpdfstring{\(2|2\)}{2|2}}
Let
\[
V=\sum_{i=1}^2 f_i(z,\xi)\frac{\partial}{\partial z_i}
 +\sum_{i=1}^2 g_i(z)\frac{\partial}{\partial\xi_i},
\qquad
\det\left(\frac{\partial g_i}{\partial z_j}(0)\right)\neq0.
\]
Write
$
B=d\bar g_1\wedge dz_1+d\bar g_2\wedge dz_2,
$ and $
C=d\bar g_1\wedge d\xi_1+d\bar g_2\wedge d\xi_2.
$
Then
\[
B^2=2d\bar g_1\wedge dz_1\wedge d\bar g_2\wedge dz_2
=-2\det(\overline{J_g})\,dz_1\wedge dz_2\wedge d\bar z_1\wedge d\bar z_2,
\]
in agreement with the sign \((-1)^{n(n+1)/2}\) for \(n=2\).
For
\[
\eta_0=\Phi\,\delta(d\xi_1)\delta(d\xi_2)
\delta(d\bar\xi_1)\delta(d\bar\xi_2),
\]
the term \(B^2/(2t^2)\) supplies the ordinary top degree, and the local residue is
\[
\Res_0(V,\eta_0)=
\left(\frac{2\pi}{i}\right)^2
\frac{[\xi_1\xi_2\bar\xi_1\bar\xi_2]\Phi(0,0,\xi,\bar\xi)}{\det J_g(0)}.
\]

For a one-lowering representative, for instance
\[
\eta_1=\Psi_1\,dz_1\,
\delta'(d\xi_1)\delta(d\xi_2)
\delta(d\bar\xi_1)\delta(d\bar\xi_2),
\]
the relevant part of \(BC\) is
\[
(d\bar g_2\wedge dz_2)(d\bar g_1\wedge d\xi_1).
\]
The factor \(d\xi_1\) lowers \(\delta'(d\xi_1)\) to \(-\delta(d\xi_1)\), while the ordinary differentials combine with the already present \(dz_1\) to give the full determinant term.  The same mechanism applies to the analogous representative with \(dz_2\delta'(d\xi_2)\).  Terms with mismatched indices, such as \(dz_1\delta'(d\xi_2)\), do not contribute through a single lowering step, as explained in Remark~\ref{rem:non-aligned-lowering}.

\section{Applications}\label{sec:applications}

We conclude with two consequences of the preceding formulae.  The first is the direct holomorphic-supergeometric analogue of the Duistermaat--Heckman principle: the integral of a \(D_V\)-closed integral form is expressed as a sum of local contributions.  The second is an elementary projective consequence of the reduced vertical component of an odd vector field on split projective superspace.  Both statements are formulated within the hypotheses established above; in particular, no distinguished geometric density is presumed to be \(D_V\)-closed without verification.

\subsection{A Duistermaat--Heckman type consequence}

Classical localisation formulae of Duistermaat--Heckman, Atiyah--Bott and Berline--Vergne express the integral of an equivariantly closed form as a sum of contributions from fixed points or fixed components \cite{DH,AB,BV}.  In the present holomorphic supergeometric setting the role of the equivariant localisation operator is played by
\[
D_V=\ddbar+i_V.
\]
Thus the corresponding statement is naturally formulated for integral forms of maximal picture which are annihilated by this operator.

\begin{corollary}\label{cor:dh-type}
Let \(S\) be a compact complex supermanifold of dimension \(n|n\), and let
\[
V\in H^0(S,T_S)_{\bar 1}
\]
be an odd holomorphic vector field whose reduced zero locus is finite and whose reduced zeros are triangular non-degenerate.  If \(\eta\) is a smooth integral form of maximal picture satisfying
\[
D_V\eta=0,
\]
then
\[
\int_S\eta
=
\sum_{p\in Z(V)_{\red}}\Res_p(V,\eta).
\]
In particular, if \(\int_S\eta\neq0\) for some \(D_V\)-closed integral form \(\eta\), then \(V\) has a reduced zero.
\end{corollary}

\begin{proof}
The equality is exactly the residue decomposition of Theorem~\ref{thm:residue-decomposition}.  If the reduced zero locus were empty, Theorem~\ref{thm:vanishing} would give \(\int_S\eta=0\).  The final assertion follows by contraposition.
\end{proof}

\subsection{Projective superspace}

Let \(\PP^{n|n}\) denote the split projective superspace in the standard sense of complex projective supergeometry \cite{Manin,Varadarajan}, whose reduced space is \(\PP^n\) and whose structure sheaf is
\[
\cO_{\PP^{n|n}}
\simeq
\bigwedge^\bullet\bigl(\cO_{\PP^n}(-1)^{\oplus n}\bigr).
\]
Equivalently, in the notation of Section~\ref{sec:conventions}, the odd conormal bundle is
\[
\cE^\vee=\cO_{\PP^n}(-1)^{\oplus n},
\qquad
\cE=\cO_{\PP^n}(1)^{\oplus n}.
\]
For an odd holomorphic vector field \(V\) on \(\PP^{n|n}\), its reduced vertical component is the section
\[
s_V\in H^0(\PP^n,\cE)
=H^0\bigl(\PP^n,\cO_{\PP^n}(1)^{\oplus n}\bigr)
\]
defined in Section~\ref{sec:conventions}.  The reduced zero locus of \(V\) is the zero scheme of this section.  The following elementary observation is the projective counterpart of the vanishing theorem.

\begin{proposition}\label{prop:projective-superspace-zero}
Every odd holomorphic vector field on the split projective superspace \(\PP^{n|n}\) has a non-empty reduced zero locus.
\end{proposition}

\begin{proof}
The section \(s_V\) is an \(n\)-tuple of global sections of \(\cO_{\PP^n}(1)\).  Hence it is represented by \(n\) homogeneous linear forms
\[
\ell_1,\ldots,\ell_n\in H^0(\PP^n,\cO_{\PP^n}(1)).
\]
Their common zero locus is the projectivisation of the kernel of the linear map
\[
\C^{n+1}\longrightarrow \C^n,
\qquad
x\longmapsto (\ell_1(x),\ldots,
\ell_n(x)).
\]
This kernel has dimension at least one.  Its projectivisation is therefore non-empty, and it is precisely the zero locus of \(s_V\).  Consequently \(Z(V)_{\red}\neq\varnothing\).
\end{proof}

\begin{corollary}\label{cor:projective-vanishing-obstruction}
No odd holomorphic vector field on \(\PP^{n|n}\) has empty reduced zero locus.
\end{corollary}

\begin{proof}
This follows immediately from Proposition~\ref{prop:projective-superspace-zero}.
\end{proof}

\begin{remark}
Proposition~\ref{prop:projective-superspace-zero} uses only the reduced vertical component of the odd vector field.  It does not require the construction of a \(D_V\)-closed projective density.  When the reduced zeros are isolated and triangular non-degenerate, the local residue formula of Theorem~\ref{thm:general-local-residue} applies to every \(D_V\)-closed integral form of maximal picture.
\end{remark}

\section*{Acknowledgements}

The first author was supported by CAPES and CONACYT and thanks CIMAT for its hospitality.  The second author was supported by CNPq grants 202374/2018-1, 302075/2015-1 and 400821/2016-8, and thanks the University of Oxford for its hospitality; he is also partially supported by the Universit\`a degli Studi di Bari  is a member of INdAM-GNSAGA

\end{document}